\documentclass[a10paper,english,17pt]{article}
\usepackage[T1]{fontenc}
\usepackage{babel}
\usepackage{graphicx,color}
\usepackage{babel}
\usepackage{amsfonts}
\usepackage{amssymb, amsmath, amsfonts, amsthm, mathrsfs}
\usepackage{hhline, latexsym}
\usepackage{float}
\usepackage[margin=3cm]{geometry}
\begin{document}
\author{Aymen Braghtha}
\title{\large \textbf{ON THE NUMBER OF ISOLATED ZEROS OF PSEUDO-ABELIAN INTEGRALS: DEGENERACIES OF THE CUSPIDAL TYPE }}
\maketitle
\begin{abstract}
We consider a multivalued function of the form $H_{\varepsilon}=P_{\varepsilon}^{\alpha_0}\prod^{k}_{i=1}P_i^{\alpha_i}, P_i\in\mathbb{R}[x,y], \alpha_i\in\mathbb{R}^{\ast}_+$, which is a Darboux first integral of polynomial one-form $\omega=M_{\varepsilon}\frac{dH_{\varepsilon}}{H_{\varepsilon}}=0, M_{\varepsilon}=P_{\varepsilon}\prod^{k}_{i=1}P_i$. We assume, for $\varepsilon=0$, that the polycyle $\{H_0=H=0\}$ has only cuspidal singularity which we assume at the origin and other singularities are saddles. 

We consider families of Darboux first integrals unfolding $H_{\varepsilon}$ (and its cuspidal point) and pseudo-Abelian integrals associated to these unfolding. Under some conditions we show the existence of uniform local bound for the number of zeros of these pseudo-Abelian integrals.
\end{abstract}
-----------------------------------------------------------------------------------------------------------------------------------
\section{Formulation of main results}
In this paper, we study a non generic case. Other non generic cases have been studied in [1,3,4,5]. Pseudo-Abelian integrals appear as the linear principal part of the displacement function in polynomial perturbation of Darboux integrable case.

More precisely consider Darboux integrable system $\omega$ given by 
\begin{equation}
\omega=M d\log H,
\end{equation}
where 
\begin{equation}
M=\prod^{k}_{i=0}P_i,\quad H=\prod^{k}_{i=0}P_i^{\alpha_i},\quad \alpha_i>0,\quad P_i\in\mathbb{R}[x,y].
\end{equation}
 
Now we consider an unfolding $\omega_{\varepsilon}$ of Darboux integrable system $\omega$, where $\omega_{\varepsilon}$ are one-forms with first integral
\begin{equation}
H_{\varepsilon}=P_{\varepsilon}^{\alpha}\prod^{k}_{i=1}P_i^{\alpha_i},\quad \omega_{\varepsilon}=M_{\varepsilon} d\log H_{\varepsilon},\quad M_{\varepsilon}=P_{\varepsilon}\prod^{k}_{i=1}P_i.
\label{moneq}
\end{equation}
where the polynomial $P_0$ has a cuspidal singularity at $p_0=(0,0)$, i.e. $P_0(x,y)=y^2-x^3+\mathcal{O}((x,y)^4)$. For non zero $\varepsilon$, the polynomial $P_{\varepsilon}=y^2-x^3-\varepsilon x^2+\mathcal{O}((x,y,\varepsilon)^4)$. 

Choose a limit periodic set i.e. bounded component of $\mathbb{R}^2\setminus\{\prod^{k}_{i=0}P_i=0\}$ filled cycles $\gamma(h)\subset\{H=h\}, h\in(0,a)$. Denote by $D\subset H^{-1}(0)$ the polycycle which is in the boundary of this limit periodic set.

Consider the unfolding $\omega_{\varepsilon}=M_{\varepsilon} d\log H_{\varepsilon}$ of the form $\omega$. The foliation $\omega_{\varepsilon}$ has a maximal nest of cycles $\gamma(\varepsilon,h)\subset\{H_{\varepsilon}=h\}, h\in(0,a(\varepsilon))$, filling a connected component of $\mathbb{R}^2\setminus\{H_{\varepsilon}=0\}$ whose boundary is a polycycle $D_{\varepsilon}$ close to $D$. Assume moreover that the foliation $\omega_{\varepsilon}=0$ has no singularities on \text{Int}$D_{\varepsilon}$.

Consider pseudo-Abelian integrals of the form 
\begin{equation}
I(\varepsilon,h):=\int_{\gamma(\varepsilon,h)}\eta_2,\quad\eta_2=\frac{\eta_1}{M_{\varepsilon}}
\label{moneq}
\end{equation}
where $\eta_1$ is a polynomial one-form of degree at most $n$.

This integral appears as the linear term with respect to $\beta$ of the displacement function of a polynomial perturbation
\begin{equation}
\omega_{\varepsilon,\beta}=\omega_{\varepsilon}+\beta\eta_1=0.
\label{moneq}
\end{equation}

We assume the following genericity assumptions 
\begin{enumerate}
\item The level curves $P_i=0, i=1,\ldots,k$ are smooth and $P_i(0,0)\neq0$.
\item The level curves $P_{\varepsilon}=0, P_i=0, i=1,\ldots,k$, intersect transversaly two by two. 
\end{enumerate}
\textbf{Theorem 1.} \emph{Under the genericity assumptions there exists a bound for the number of isolated zeros of the pseudo-Abelian integrals $I(\varepsilon,h)=\int_{\gamma(\varepsilon,h)}\eta_2$ in $(0,a(\varepsilon))$. The bound is locally uniform with respect to all parameters in particular in $\varepsilon$.}\\
\linebreak
Let $\mathcal{F}_1:\{\omega_{\varepsilon}=0\}, \mathcal{F}_2:\{d\varepsilon=0\}$ are the foliations of dimension two in complex space of dimension three with coordinates $(x,y,\varepsilon)$.

Let $\mathcal{F}$ be the foliation of dimension one on the complex space of dimension three with coordinates $(x,y,\varepsilon)$ which is given by the intersection of leaves of $\mathcal{F}_1$ and $\mathcal{F}_2$ (i.e. given by the 2-form $\Omega=\omega_{\varepsilon}\wedge d\varepsilon$). This foliation has a cuspidal singularity at the origin (a cusp).

We want to study the analytical properties of the foliation $\mathcal{F}$ in a neighborhood of the cusp. For this reason we make a global blowing-up of the cusp of the product space $(x,y,\varepsilon)$ of phase and parameter spaces. We want our blow-up to seperate the two branche of the cusp. This requirements leads to the quasi-homogeneous blowing-up of weight $(2,3,2)$.\\
\linebreak
\textbf{Remark 1.} \emph{In term of first integrals, the foliation $\mathcal{F}$ is given by two first integrals 
$$H(x,y,\varepsilon)=h,\quad\varepsilon=s.$$} 

\section{Quasi-homogeneous blowing-up of $\mathcal{F}$}

Recall the construction of the quasi-homogenous blowing-up. We define the weighted projective space $\mathbb{CP}^2_{2:3:2}$ as factor space of $\mathbb{C}^3$ by the $\mathbb{C}^{\ast}$ action $(x,y,\varepsilon)\mapsto(t^2x,t^3y,t^2\varepsilon)$. The quasi-homogeneous blowing-up of $\mathbb{C}^3$ at the origin is defined as the incidence three dimensional manifold $W=\{(p,q)\in\mathbb{CP}^{2}_{2:3:2}\times\mathbb{C}^{3}:\exists t\in\mathbb{C}: (q_1,q_2,q_3)=(t^2p_1,t^{3}p_2,t^{2}p_3)\}$, where $(q_1,q_2,q_3)\in\mathbb{C}$ and $[(p_1,p_2,p_3)]\in\mathbb{CP}^{2}_{2:3:2}$.

The quasi-homogeneous blowing-up $\sigma: W\rightarrow\mathbb{C}^3$ is just the restriction to $W$ of the projection $\mathbb{CP}^{2}_{2:3:2}\times\mathbb{C}^{3}$.

We will need explicit formula for the blow-up in the standard affine charts of $W$. The projective space $\mathbb{CP}^2_{2:3:2}$ is covred by three affine charts: $U_1=\{x\neq0\}$ with coordinates $(y_1,z_1)$, $U_2=\{y\neq0\}$ with coordintaes $(x_2,z_2)$ and $U_3=\{\varepsilon\neq0\}$ with coordinates $(x_3,y_3)$. 

The transition formula follow from the requirement that the points $(1,y_1,z_1), (x_2,1,z_2)$ and $(x_3,y_3,1)$ lie on the same orbit of the action:
\begin{align*}
&F_2:(y_1,z_1)\mapsto\left(x_2=1/y_1^{2/3}, z_2=z_1/y_1\sqrt{y_1}\right)\\
&F_3:(y_1,z_1)\mapsto\left(x_3=1/{z_1}, y_3=y_1/z_1\sqrt{z_1}\right).
\end{align*}

These affine charts define affine charts on $W$, with coordinates $(y_1,z_1,t_1),\\ (x_2,z_2,t_2)$ and
 $(x_3, y_3,t_3)$. The blow-up $\sigma$ is written as
\begin{align}
&\sigma_1:\quad x=t_1^2,\qquad y=t_1^3y_1,\qquad \varepsilon=t_1^2z_1\\
&\sigma_2:\quad x=t_2^2x_2,\qquad y=t_2^3,\qquad \varepsilon=t^2_2z_2\\
&\sigma_3:\quad x=t^2_3x_3,\qquad y=t_3^3y_3,\qquad \varepsilon=t_3^2.
\end{align}

We apply this blow-up $\sigma$ to the one-dimensional foliation $\mathcal{F}$. Let $\sigma^{-1}\mathcal{F}$ the lifting of the foliation $\mathcal{F}$ to the complement  This foliation has a cuspidal singularity at the origin. The pull-back foliation $\sigma^{\ast}\mathcal{F}$ will be called the strict transform of the foliation $\mathcal{F}$ is defined by the pull-back $\sigma^{\ast}\Omega=\sigma^{\ast}$ ($\omega_{\varepsilon}\wedge d\varepsilon$) divided by a suitable power of the function defining the exceptional divisor. In this charts $U_j, j=1,2,3$ we have 
$$
\sigma^{\ast}_1\Omega=x^2\Omega_1,\quad
\sigma^{\ast}_2\Omega=y^3\Omega_2,\quad
\sigma^{\ast}_3\Omega=\varepsilon^2\Omega_3,
$$
where 
\begin{align}
&\Omega_{1}=(6y_1^2-6-4z_1)dx\wedge dz_1+4y_1z_1dy_1\wedge dx+2xy_1dy_1\wedge dz_1,\\
&\Omega_2=(6-6x_2^3-4x^2_2z_2)dy\wedge dz_2+(-6z_2x_2^2-4x_2z_2^2)dx_2\wedge dy\\
&+(-3yx_2^2-2yx_2z_2)dx_2\wedge dz_2,\\
&\Omega_3=(-6x_3^2-4x_3)dx_3\wedge d\varepsilon+4y_3dy_3\wedge d\varepsilon.
\end{align}
\textbf{Remark 2.} \emph{In term of first integrals, the foliation $\sigma^{\ast}\mathcal{F}$ is given by two first integrals 
$$
\sigma^{\ast}H(x,y,\varepsilon)=h,\quad\sigma^{\ast}\varepsilon=s,
$$
In particular in a neighborhood of the exceptional divisor the restrictions of the foliation $\sigma^{\ast}\mathcal{F}$ to the charts $U_1$ and $U_3$ are given respectively, by
\begin{align}
&\psi_1=H(t_1^2,t_1^3y_1,t_1^2z_1)=x^3(y_1^2-1)=h,\quad\varphi_1=xz_1=s,\\
&\psi_3=H(t_3^2x_3,t_3^3y_3,t_3^2)=\varepsilon^3(y_3^2-x_3^2-x_3^3)=h,\quad\varphi_3=\varepsilon=s,
\end{align}
where $\{x=0\}$ and $\{\varepsilon=0\}$ are local equations of the exceptional divisor respectively.}

\section{Singular locus of the foliation $\sigma^{\ast}\mathcal{F}$}

In this section, we compute the singular locus of the pull-back $\sigma^{\ast}\Omega$ in a neighborhood of the exceptional divisor $\mathbb{CP}^{2}_{2:3:2}$. We check it in each chart seperatly.

In the chart $U_1$, the zeros locus of the form $\Omega_1$ in a neighborhood of the exceptional divisor $\{x=0\}$ consists of germs of two curves $\{y_{1}=\pm1,z_1=0\}$ and a two singular points $p_1=(0,1,0), p_2=(0,-1,0)$ generated by the quasi-homogeneous blowing-up. 

In the chart $U_3$, the zeros locus of the form $\Omega_3$ in a neighborhood of the exceptional divisor $\{\varepsilon=0\}$ consists of $p_3=(0,0,0)$ (Morse point) and $p_4=(-\frac{2}{3},0,0)$ (center). The singularities of this foliation are the line of Morse points $x_3=0, y_3=0$, the lines of centers $x_3=-\frac{2}{3},y_3=0$ and the transform strict of $\{y^2-x^3-x^2\varepsilon=0\}$.\\
\linebreak
\textbf{Proposition 1.} \emph{The singularities of $\sigma^{\ast}\mathcal{F}$ are located at the points $ p_1, p_2,p_3$ and $p_4$. The points $p_1,p_2$ and $p_3$ are linearisable saddles and the point $p_4$ is a center.}
\begin{proof} 
Since $\sigma:W\rightarrow\mathbb{C}^3$ is a biholomorphism autside the exceptional divisor $\mathbb{CP}^2_{2:3:2}$, all singularities of $\sigma^{\ast}\mathcal{F}$ on $\mathbb{C}^3\setminus\{x=0\}$ correspond to singularities of $\mathcal{F}$. Thus, it suffices to compute the singularities of $\sigma^{\ast}\mathcal{F}$ on the exceptional divisor $\{x=0\}$. More precisely, we consider the foliation on neighborhood of $\mathbb{CP}^2_{(2:3:2)}$ (the exceptional divisor) generated by the blown-up one-form $\sigma^{\ast}\Omega$. Let $\psi_1,\psi_3$ are the functions given in (13) and (14).\\
\linebreak
(1) In the chart $U_1$, near the divisor exceptional and for $|z_1|\leq\epsilon$ for $\epsilon$ sufficiently small, the foliation $\sigma^{\ast}\mathcal{F}$ is given by two first integrals
$$G_1=\varphi_1^3\psi_1^{-1}=z_1^3(y_1^2-(1+z_1))^{-1}V^{-1}=s^3h^{-1},\quad\varphi_1=xz_1=s.$$ 
where $V$ is analytic function such that $V(0,0,0)\ne0$. In particular on the exceptional divisor $\{x=0\}$ the foliation $\sigma^{\ast}\mathcal{F}$ is given by the levels $G_1=s^3h^{-1}=t$.\\ 
\linebreak
Now we calculate the eigenvalues at $p_1$ and $p_2$. The vector field $V_1$ generating the foliation $\sigma^{\ast}\mathcal{F}$ is given by
$$
V_1(x,y_1,z_1)=\beta_1x\frac{\partial}{\partial x}+\beta_2y_1\frac{\partial}{\partial y_1}+\beta_3z_1\frac{\partial}{\partial z_1},
$$
where the vector $(\beta_1,\beta_2,\beta_3)$ satisfies the following equations
$$
<(\beta_1,\beta_2,\beta_3),(3,1,0)>=0, <(\beta_1,\beta_2,\beta_3),(1,0,1)>=0
$$
here $<,>$ be the usual scalar product on $\mathbb{C}^3$. By simple computation, we obtain $\beta_1=1,\beta_2=-3$ and $\beta_3=-1$.\\
\linebreak
(2) In the chart $U_3$, near the exceptional divisor $\{\varepsilon=0\}$, the foliation $\sigma^{\ast}\mathcal{F}$ is given by 
$$
G_3=\varphi_3^3\psi_3^{-1}=(y_3^2-x_3^2(1+x_3))^{-1}=s^3h^{-1},\quad\varphi_3=\varepsilon=s.
$$
In particular the restriction of this foliation to the exceptional divisor $\{\varepsilon=0\}$, by Morse lemma we can put the function $1/G_3$ to the normal form $y_3^2-z^2_3$ in a neighborhood of $p_3$ (we put the variable change $z_3=\pm x_3(1+x_3)^{1/2}$). On other hand the Hessian matrix of  $1/G_3$ at the point $p_4$ has two positive eigenvalus.   
\end{proof}

\section{The different scaled variations of $\delta(s,t)$}

In this section, we compute the scaled variations with respect to differents variables $s$ and $t$ of the integrals of the blown- up one form $\sigma^{\ast}_1\eta_2$ along the different relatives cycles using the same technics of [5].\\
\linebreak
\textbf{Proposition 2.} \emph{The computation of the different scaled variations of the cycle $\delta(s,t)$ us gives
\begin{enumerate}
\item For $t\in[0,2N]$, the cycle $\delta(s,t)$ satisfies a iterated scaled variations with respect to $t$ of the form
\begin{equation}
\mathcal{V}ar_{(t,3)}\circ\mathcal{V}ar_{(t,-1)}\circ\mathcal{V}ar_{(t,-\alpha_1)}\circ\ldots\circ \mathcal{V}ar_{(t,-\alpha_k)}\delta(s,t)=0.
\end{equation} 
\item For $t\in[N,+\infty]$, the cycle $\delta(s,t)$ satisfies a iterated scaled variations with respect to $1/t$ of the form
\begin{equation}
\mathcal{V}ar_{(1/t,-3)}\circ \mathcal{V}ar_{(1/t,1)}\circ \mathcal{V}ar_{(1/t,1)}\circ\mathcal{V}ar_{(1/t,\alpha_1)}\circ\ldots\circ\mathcal{V}ar_{(1/t,\alpha_k)}\delta(s,1/t)=0.
\end{equation}
\item Near $s=0$, we have
\begin{equation}
\mathcal{V}ar_{(s,1)}\circ\mathcal{V}ar_{(s,1)}\delta(s,t)=\mathcal{V}ar_{(s,1)}(\tilde{\delta}(s,t))=0,
\end{equation}
where $\mathcal{V}ar_{(s,1)}\delta(s,t)=\tilde{\delta}(s,t)$ is a figure eight cycle.
\end{enumerate}}
\begin{proof} As in [5], there exist a some local chart with coordinates $(u,v,w)$ defined in a some neighborhood of each separatrix of polycycle such that the foliation $\sigma^{\ast}\mathcal{F}$ is defined by two first integrals. Precisely:
\begin{enumerate}
\item for $t\in[0,2N]$, there exist a local chart $(V_{div},(u,v,w))$ defined in neighborhood of the separatrix $\delta_{div}$ such the foliation $\sigma^{\ast}\mathcal{F}$ by two first integrals
$$
F_1=w^3(v-1)^{-1}(v+1)^{-1}=t,\quad F_2=uw=s,
$$
\item for $t\in[N,+\infty]$, there exists a local chart $(V^{+}_{div},(u,v,w))$ defined in neighborhood of the separatrix $\delta^{+}_{div}$ such that the foliation $\sigma^{\ast}_1\mathcal{F}$ is defined by two first integrals
$$
F_1=w^3(v+2)^{-1}v^{-1}=t,\quad F_2=uw=s,
$$
\item for $t\in[N,+\infty]$, there exists a local chart $(V^{-}_{div},(u,v,w))$ defined in a neighborhood of the separatrix $\delta^{-}_{div}$ such that the foliation $\sigma^{\ast}_1\mathcal{F}$ is defined by two first integrals
$$
F_1=w^3(v-2)^{-1}v^{-1}=t,\quad F_2=uw=s.
$$
\end{enumerate}
 In second step we prove that each relative cycle can be chosen as a lift of a path contained in the separatrix associated to this relative cycle. Precisely:
\begin{enumerate}
\item on the chart $(V_{div},(u,v,w))$, the linear projection $\pi(u,v,w)=v$ is every where transverse to the levels of the foliation $\sigma^{\ast}\mathcal{F}$ which corresponds simply to the graphs of the multivalued functions
 $$v\mapsto (u,w)=\left(st^{-\frac{1}{3}}(v-1)^{-\frac{1}{3}}(v+1)^{-\frac{1}{3}},t^{\frac{1}{3}}(v-1)^{\frac{1}{3}}(v+1)^{\frac{1}{3}}\right),$$
\item on the chart $(V^{+}_{div},(u,v,w))$, the linear projection $\pi(u,v,w)=v$ is every where transverse to the levels of the foliation $\sigma^{\ast}\mathcal{F}$ which corresponds simply to the graphs of the multivalued functions
 $$v\mapsto (u,w)=\left(st^{-\frac{1}{3}}v^{-\frac{1}{3}}(v+2)^{-\frac{1}{3}},t^{\frac{1}{3}}v^{\frac{1}{3}}(v+2)^{\frac{1}{3}}\right),$$
\item on the chart $(V^{-}_{div},(u,v,w))$, the linear projection $\pi(u,v,w)=v$ is every where transverse to the levels of the foliation $\sigma^{\ast}\mathcal{F}$ which corresponds simply to the graphs of the multivalued functions
 $$v\mapsto (u,w)=\left(st^{-\frac{1}{3}}v^{-\frac{1}{3}}(v-2)^{-\frac{1}{3}},t^{\frac{1}{3}}v^{\frac{1}{3}}(v-2)^{\frac{1}{3}}\right).$$
\end{enumerate}
In third step, we compute the different scaled variations of relatives cycles using the local expression of two first integrals $F_1$ and $F_2$ above near the singular points $p_1, p_2$ and $p_3$. Recall that the scaled variation of a relative cycle $\delta(s)$ is given by
$$
\mathcal{V}ar_{(s,\beta)}\delta(s)=\delta(se^{i\pi\beta})-\delta(se^{-i\pi\beta}).
$$ 
In the local chart $(V^+_{div},(u,v,w))$, the restriction of the blown-up foliation $\sigma^{\ast}_1\mathcal{F}$ to the transversals sections $\Sigma^{-}_{div}=\{w=1\}$ (near the point $p_3$) and $\Omega_{+}=\{u=1\}$ (near the point $p_1$) is given respectively by  
\begin{align*}
&F_1|_{\Sigma_{div}^-}=\frac{1}{v}=t,\quad F_2|_{\Sigma_{div}^-}=u=s,\\
&F_1|_{\Omega_{+}}=\frac{w^3}{v}=t,\quad F_2|_{\Omega_{+}}=w=s.
\end{align*}

Let us fix $t\in[N,+\infty]$. We observe that the restriction of the foliation $\sigma^{\ast}_1\mathcal{F}$ to the transversal section $\Sigma^{+}_{div}=\{w=1\}$ is analytic with respect to $s$. Then, after taking an scaled variation with respect to $s$, the relative cycle $\delta^{+}_{div}(s,t)$ is replaced by a loop $\theta_1$, modulo homotopy, which consists of line segment $\ell_{31}=[p_3,p_1]$ connecting the Morse point $p_3$ with the point $p_1$ encircling the latter along a small counterclockwise circular arc $\alpha_1$ and then returning along the segment  $\ell_{13}=[p_1,p_3]$. The loop $\theta_1$ can be moved along the complex curve $\{u=w=0\}$. Then, we have 
$$
\mathcal{V}ar_{(s,1)}\delta^{+}_{div}(s,t)=\theta_1=\ell_{31}\alpha_{1}\ell_{13}.
$$
The same computation of the scaled variation with respect to $s$ for the relative cycle  $\delta^{-}_{div}(s,t)$ gives us a loop $\theta_3$, modulo homtopy, which can be moved along the complex plane $\{u=w=0\}$. The loop $\theta_3$ consists of line segment $\ell_{32}=[p_3,p_2]$ connecting the point $p_3$ with the point $p_2$ encircling the latter along a small counterclockwise circular arc $\alpha_3$ and then returning along the segment  $\ell_{23}=[p_2,p_3]$. Then, we have
$$
\mathcal{V}ar_{(s,1)}\delta^{-}_{div}(s,t)=\theta_3=\ell_{32}\alpha_{3}\ell_{23}.
$$

In the local chart $(V_{div},(u,v,w))$, we define the transversal section $\Omega_{+}=\{u=1\}$ (resp $\Omega_{+}=\{u=1\}$) near $p_1$ (resp near $p_2$). The restriction of the foliation $\sigma^{\ast}_1\mathcal{F}$ to the transversal section $\Omega_+$ is given by
$$
F_1|_{\Omega_+}=\frac{w^3}{v}=t,\quad F_2|_{\Omega_+}=w=s.
$$ 
On the second step let us fix $t\in[0,2N]$. After taking an scaled variation with respect to $s$, the relative cycle $\delta_{div}(s,t)$ is replaced by a figure eight cycle which can be moved along the complex line $C^t_{div}=\{x=0,G_1=t\}$ of the foliation $\sigma^{\ast}_1\mathcal{F}$. This case is similar to the classical situation which is studied by Bobie\'{n}ski and Marde\v{s}i\'{c} in [2].

Now using the analycity of the lifting $\sigma^{-1}\mathcal{F}$ with respect to $s$, the scaled variation of the cycle of integration $\delta(s,t)$ with respect to $s$ is equal to the scaled variation with respect to $s$ of the following difference $\delta_{div}^{+}(s,t)-\delta_{div}^{-}(s,t)$ which is equal, modulo homotopy, to the cycle $\theta_1\theta_3^{-1}$, where $\theta_3^{-1}$ is the inverse of the loop $\theta_3$. Shematically, the loop $\theta_1\theta_3^{-1}$ is a figure eight cycle.
\end{proof}
\textbf{Remark 3.} 
\begin{itemize}
\item In the local chart $(V^+_{div},(u,v,w))$ (resp $V^-_{div},(u,v,w)$), the loop $\theta_1$ (resp $\theta_3$) generating the fundamental group of the complex plane  $\{u=w=0\}\setminus\{p_1\}$ (resp $\{u=w=0\}\setminus\{p_2\}$) with base point $p_3$.
\item  By the univalness of the blown-up one form $\sigma^{\ast}_1\eta_2$, we have 
$$
\mathcal{V}ar_{(t,\alpha)}\int_{\delta(s,t)}\sigma^{\ast}_1\eta_2=\int_{\mathcal{V}ar_{(t,\alpha)}\delta(s,t)}\sigma^{\ast}_1\eta_2.
$$
\end{itemize}


\section{Proof of the Theorem}

In this section we first take benefit from the blowing-up in the family to prove our principal theorem. the proof is analoguous of the following :\\
\linebreak
 \textbf{Theorem 2.} \emph{There exists a bound of the number of zeros of the function $t\mapsto J(s,t)$, for $t\in[0,+\infty]$ and $s>0$ sufficiently small. This bound is locally with respect to all parameters uniform, in particular with respect to $s$.}\\

Let $\beta=(\beta_1,\dots,\beta_{k+2})$ where $\beta_1=3, \beta_2=-1,\beta_3=-\alpha_1,\ldots,\beta_{k+2}=-\alpha_k$. Let $D_1$ is slit annulus in the complex plane $\mathbb{C}^{\ast}_t$ with boundary $\partial D_1$. This boundary is decomposed as follows $\partial D_1= C_{R_1}\cup C_{r_1}\cup C^{\pm}$, where 
$C_{R_1}=\{|t|=R_1, |\arg t|\leq \alpha\pi\}, C^{\pm}=\{r_1<|t|<R_1, |\arg t|=\pm\alpha\}$ and $C_{r_1}=\{|t|=r_1, |\arg t|\leq\alpha\pi\}$.\\

Petrov's method gives us that the number of zeros $\#Z(J(s,t))$ of the function $J(s,t)$ in slit annulus $D_1$ is bounded by the increment of the argument of $J(s,t)$ along $\partial D_1$ divided by $2\pi$ i.e.
\begin{align*}
\#Z(J(s,t)|_{D_1})\leq&\frac{1}{2\pi}\Delta\arg(J(s,t)|_{\partial D_1})=\frac{1}{2\pi}\Delta\arg(J(s,t)|_{C_{R_1}})\\
&+\frac{1}{2\pi}\Delta\arg(J(s,t)|_{C^{\pm}})+\frac{1}{2\pi}\Delta\arg(J(s,t)|_{C_{r_1}})
\end{align*}
\linebreak
(\textbf{A}) The increment of argument $\Delta\arg(J(s,t)|_{C_{R_1}})$ is uniformly bounded by Gabrielov's theorem [6]. \\
\linebreak
(\textbf{B}) We use the Schwartz's principle
$$
Im(J(s,t)) |_{C^{\pm}}=\mp2i\mathcal{V}ar_{(t,\alpha)}J(s,t).
$$ 
Thus, the increments of argument along segments $C^{\pm}$ are bounded by zeros of the variation $\mathcal{V}ar_{(t,\alpha)}J(s,t)$ on segment $(r,R)$. By identity (18), the function $\mathcal{V}ar_{(t,\beta_i)}J(s,t)$ can be written as follows 
\begin{align*}
\mathcal{V}ar_{(t,\beta_i)}J(s,t)&=K(t^{\frac{\beta_1}{\beta_i}},\dots,t^{\frac{\beta_{k+\mu}}{\beta_i}},s;\log s)\\
&=K(e^{\frac{\beta_1}{\beta_i}\log t},\dots,e^{\frac{\beta_{k+\mu}}{\beta_i}\log t},e^{\log s};\log s)
\end{align*} 
where $K$ is a meromorphic function. The function $\mathcal{V}ar_{(t,\beta_i)}J(s,t)$ is logarithmico-analytic function of type 1 in the variable $s$ (see [9]). Then, there exist a finit recover of $\mathbb{R}^{k+\mu+1}\times\mathbb{R}$ by a logarithmico-exponential cylinders, using Rolin-Lion's theorem [9], such that on each cylinder of this family we have
$$
\mathcal{V}ar_{(t,\beta_i)}J(s,t)=y_0^{r_0}y_1^{r_1}A(t)U(t,y_0,y_1),
$$ 
with $y_0=s-\theta_0(t), y_1=\log y_0 -\theta_1(t)$, where $\theta_0, \theta_1, A$ are logarithmico-exponential functions and $U$ is a logarithmico-exponential unity function. As the number of zeros of a logarithmico-exponential function is bounded, the number of zeros of $\mathcal{V}ar_{(t,\beta_i)}J(s,t)$ is bounded. \\
\linebreak
(\textbf{C}) Finally, we estimate the increment of argument of $J$ along the small arc $C_{r_1}$. Then, it is necessarily to study the increment of argument of the leading term of the function $J$ at $t=0$.\\
\linebreak
\textbf{Lemma 1.} \emph{The increment of the argument of $J(s,t)$ along the small circular arc $C_{r_1}$ can be estimated by the increment of the argument of a some meromorphic function  $F(s,t)$.}
\begin{proof}
The problem of the estimation of the increment of the argument of $J(s,t)$ along the circular arc $C_{r_1}$ consist that the principal part of the function $J$ contains the term $\log s\rightarrow-\infty$ as $s\rightarrow0$. To resolve this problem we make a blowing-up at the origin in the total space with coordinates $(x,y,z)$ where 
 $$
x=J_1(s,t),\quad y=J_2(s,t),\quad z=(\log s)^{-1}. 
$$
The function $J(s,t)$ can be rewritten as follows
\begin{equation*}
J(s,t)=J_1(s,t)+J_2(s,t)\log s=((\log s)^{-1}J_1(s,t)+J_2(s,t))\log s=(z x+y)z^{-1}.
\end{equation*}
Thus, for $z^{-1}\in\mathbb{R}$ be sufficiently small, we have 
$$
\arg( J(s,t))=\arg((zx+y)z^{-1})=\arg(zx+y).
$$

To estimate the increment of argument of $zx+y$ uniformly with respect to $s>0$ we make a quasi-homogeneous blowing-up $\pi_1$ with weight $(\frac{1}{2},1,\frac{1}{2})$ of the polynomial $zx+y$ at $C_1=\{x=y=z=0\}$ (the centre of blowing-up). The explicit formula of the quasi-homogeneous blowing-up $\pi_1$ in the affine charts $T_1=\{x\ne0\}, T_2=\{y\ne0\}$  and $T_3=\{z\ne0\}$ is written respectively as
\begin{align*}
&\pi_{11}: x=\sqrt{x_1},\quad y=y_1x_1,\quad z=z_1\sqrt{x_1},\\
&\pi_{12}: x=x_2\sqrt{y_2},\quad y=y_2,\quad z=z_2\sqrt{y_2},\\
&\pi_{13}: x=x_3\sqrt{z_3},\quad y=y_3z_3,\quad z=\sqrt{z_3}.
\end{align*}
The pull-back $\pi^{\ast}_{1}(zx+y)$ is given, in different charts, by
\begin{align*}
&\pi^{\ast}_{11}(zx+y)=x_1(z_1+y_1)=d_1P_1(x_1,y_1,z_1),\\
&\pi^{\ast}_{12}(zx+y)=y_2(x_2z_2+1)=d_2P_2(x_2,y_2,z_2),\\
&\pi^{\ast}_{13}(zx+y)=z_3(x_3+y_3)=d_3P_3(x_3,y_3,z_3).
\end{align*}
where $d_i=0$ and $P_i=0$ are equations of exceptional divisor and the strict transform of $zx+y=0$ respectively.

Observe that $P_i=0, i=1,3$ has not a normal crossing with tha exceptional divisor $d_i=0, i=1,3$. To resolve this problem we make a second blowing-up $\pi_2$ with centre a subvariety $C_2$ which is given, in differents charts, as following:
\begin{enumerate}
\item In the chart $T_1$, choose a local coordinate chart with coordinates $(x_1,y_1,z_1)$ in which $C_2=\{ y_1=z_1=0\}$. Then $\pi_2^{-1}(C_2)$ is covred by two coordinates charts $U_{y_1}$ and $U_{z_1}$ with coordinate $(\tilde{x}_1,\tilde{y}_1,\tilde{z}_1)$ where in $y_{1}$-chart $U_{y_1}$ the blowing-up $\pi_2$ is given by $x_1=\tilde{x}_1, y_1=\tilde{y}_1, z_1=\tilde{z}_1\tilde{y}_1$ and in $z_{1}$-chart $U_{z_1}$ the blowing-up $\pi_2$ is given by $x_1=\tilde{x}_1, y_1=\tilde{y}_1\tilde{z}_1, z_1=\tilde{z}_1$.\\
\item In the chart $T_2$, the blowing-up $\pi_2$ is a biholomorphism ($\pi_2$ is a proper map).\\ 
\item In this chart $T_3$, choose a local coordinate chart with coordinates $(x_3,y_3,z_3)$ in which $C_2=\{x_3=y_3=0\}$. Then $\pi_2^{-1}(C_2)$ is covred by two coordinates charts $U_{x_3}$ and $U_{y_3}$ with coordinate $(\tilde{x}_3,\tilde{y}_3,\tilde{z}_3)$ where in $x_{3}$-chart $U_{x_3}$ the blowing-up $\pi_2$ is given by  $x_3=\tilde{x}_3, y_3=\tilde{y}_3\tilde{x}_3, z_3=\tilde{z}_3$ and in $y_{3}$-chart $U_{y_3}$ the blowing-up $\pi_2$ is given by $x_3=\tilde{x}_3\tilde{y}_3, y_3=\tilde{y}_3, z_3=\tilde{z}_3$.
\end{enumerate}
The pull-back $\pi_1^{\ast}(zx+y)$ is given, in different charts, by\
\begin{itemize}
\item In the $y_{1}$-chart $U_{y_1}$, the transformation of the pull-back $\pi_1^{\ast}(zx+y)$ by the blowing-up $\pi_2$ is given by  
$$\pi_{2}^{\ast}\circ\pi_1^{\ast}(zx+y)=\pi_{2}^{\ast}(d_1P_1(x_1,y_1,z_1))=\tilde{x}_1\tilde{y}_1(\tilde{z}_1+1)\overset{0}{\approx}\tilde{x}_1\tilde{y}_1=J_2(s,t)=F(s,t).$$
\item In the $z_{1}$-chart $U_{z_1}$, the transformation of the pull-back $\pi_1^{\ast}(zx+y)$ by the blowing-up $\pi_2$ is given by 
$$\pi_{2}^{\ast}\circ\pi_1^{\ast}(zx+y)=\pi_{2}^{\ast}(d_1P_1(x_1,y_1,z_1))=\tilde{z}_1\tilde{x}_1(\tilde{y}_1+1)\overset{0}{\approx}\tilde{x}_1\tilde{z}_1=(\log s)^{-1}J_1(s,t)=F(s,t).$$
\item In the chart $T_2$, we have
$$\pi_{2}^{\ast}\circ\pi_1^{\ast}(zx+y)=\pi_2^{\ast}(d_2P_2(x_2,y_2,z_2))=d_2P_2(x_2,y_2,z_2)=(\log s)^{-1}J_1(s,t)+J_2(s,t)=F(s,t).$$
\item In the $x_{3}$-chart $U_{x_3}$,  the transformation of the pull-back $\pi_1^{\ast}(zx+y)$ by the blowing-up $\pi_2$ is given by  
$$\pi_{2}^{\ast}\circ\pi_1^{\ast}(zx+y)=\pi_{2}^{\ast}(d_3P_3(x_3,y_3,z_3))=\tilde{z}_3\tilde{x}_3(\tilde{y}_3+1)\overset{0}{\approx}\tilde{x}_3\tilde{z}_3=(\log s)^{-1}J_1(s,t)=F(s,t).$$
\item In the $y_{3}$-chart $U_{y_3}$, the transformation of the pull-back $\pi_1^{\ast}(zx+y)$ by the blowing-up $\pi_2$ is given by 
$$\pi_{2}^{\ast}\circ\pi_1^{\ast}(zx+y)=\pi_{2}^{\ast}(d_3P_3(x_3,y_3,z_3))=\tilde{z}_3\tilde{y}_3(\tilde{x}_3+1)\overset{0}{\approx}\tilde{y}_3\tilde{z}_3=J_2(s,t)=F(s,t).$$
\end{itemize}
Finally, we distinguish three cases:
\begin{enumerate}
\item $\arg_{C_{r_1}}J(s,t)=\arg_{C_{r_1}}((\log s)^{-1}J_1(s,t))=\arg_{C_{r_1}}J_1(s,t)$, ($(\log s)^{-1}\in\mathbb{R}$)
\item $\arg_{C_{r_1}}J(s,t)=\arg_{C_{r_1}}J_2(s,t)$,
\item In the chart $T_2$, the function $F(s,t)=((\log s)^{-1}J_1(s,t))+J_2(s,t)$ is meromorphic.
\end{enumerate}
\end{proof}
Now we define the functional space $\mathcal{P}_{\beta}$ which are formed of coefficients of the polynomials $P_i$ of the Darboux first integral $H$, the coefficients of the polynomials $R, S$ of the perturbative one forme $\eta$, exponents $\alpha_i$ and degrees $n_i=\deg P_i, n=\max(\deg R, \deg S)$. Consider the following finite dimensional functional space $\mathcal{P}_{\beta}$
\begin{align*}
\mathcal{P}_{\beta}(m_{\beta},M_{\beta};\beta_1,\ldots,\beta_{k+2})=\{&\sum_{j=1}^{k+2}\sum_{n,\ell}A_{j\ell n}(s) t^{\beta_jn}s^{m}\log^{\ell}(t):\\
&A_{j\ell n}(s) \in\mathbb{C},
m_{\beta}<A_{j\ell n}<M_{\beta}, 0\le \ell\le k+1\}.
\end{align*}
For the first two cases, the function $J_i(s,t), i=1,2$ satisfies the following iterated variations equation with respect to $t$
$$\mathcal{V}ar_{(t,\beta_1)}\circ \ldots\circ\mathcal{V}ar_{(t,\beta_{k+2})}J_i(s,t)=0.$$
Thus, by Lemma 4.8 from [2], there exists a non zero leading term $P_{i\beta}\in\mathcal{P}_{\beta}$ of $J_i(s,t)$, i=1,2 at $t=0$ such that $|J_i(s,t)-P_{i\beta}(s,t)|=O(t^{\mu_1}), \mu_1>0$, uniformly in $s$. Moreover, the function $J_i(s,t), i=1,2$ satisfies the iterated variation equation
$$\mathcal{V}ar_{(s,1)}J_i(s,t)=0.$$
Thus, we have $J_i(s,t)=O(s^{\mu_2}), \mu_2>0$, uniformly in $t$.

For each element in the parameter space, we can choose the leading term of $P_{i\beta}$. The increment of argument of this leading term is bounded by a constant $C(M_{\beta},k+2,\beta_{k+2})$. Since the leading term of $P_{i\beta}$ is also the leading term of $J_i(s,t)$, the limit $\lim_{r_1\rightarrow 0}\Delta\arg(J_i(s,t)|C_{r_1})\leq C(M_{\beta},k+2,\beta_{k+2})$.\\

In the chart $T_2$, the function $F$ is meromorphic. Thus, this function can be rewritten as following
$$
F(s,t)=(\log s)^{-1}J_1(s,t)+J_2(s,t)=G(t^{\beta_1},\ldots,t^{\beta_k},s,(\log s)^{-1})
$$
where $G$ is meromorphic function. The number $\#Z(G)$ of zeros of the function $G$ is uniformly bounded. The latter claim is a direct application of fewnomials theory of Khovanskii [8]: since the functions $\epsilon_i(t)=t^{\beta_i}, \epsilon(s)=(\log s)^{-1}$ are Pfaffian functions (solutions of Pfaffian equations $td\epsilon_i-\beta_i \epsilon_i dt=0$ and $sd\epsilon+\epsilon^2 ds$, respectively), the upper bound for this number of zeros can be given, using Rolle-Khovanskii arguments of [7], in terms of the number of zeros of some polynomial and its derivatives. The latter are uniformly bounded by Gabrielov's theorem [6].

\textbf{Burgundy University, Burgundy Institue of Mathematics,\\
U.M.R. 5584 du C.N.R.S., B.P. 47870, 21078 Dijon \\
Cedex - France.}\\
\textbf{E-mail adress:} aymenbraghtha@yahoo.fr

\end{document}